\documentclass[11pt,titlepage]{article}
\usepackage[margin=1.2in]{geometry}
\usepackage{amsmath, amsbsy, amsthm, array, mathrsfs}
\usepackage{graphics}
\usepackage{physics}
\usepackage{epsfig, amssymb,latexsym,verbatim}
\usepackage{graphicx}
\usepackage{color}
\usepackage{dsfont, bbm}
\usepackage{relsize}

\newcommand{\R}{\mathbb{R}}
\newcommand{\mbs}{\mathbb{S}}

\newcommand{\noin}{\noindent}
\newcommand{\bee}{\begin{eqnarray*}}
\newcommand{\ene}{\end{eqnarray*}}
\newcommand{\bec}{\begin{center}}
\newcommand{\enc}{\end{center}}
\newcommand{\be}{\begin{equation}}
\newcommand{\ee}{\end{equation}}

\newcommand{\ep}{\varepsilon}
\newcommand{\mb}{\mathbf}
\newcommand{\bs}{\boldsymbol}
\newcommand{\tb}{\textbf}
\newcommand{\pend}{$\blacksquare$}
\newcommand{\vs}{\vskip 3mm}
\newcommand{\bi}{\begin{itemize}}
\newcommand{\ei}{\end{itemize}}
\begin{document}

\title{\LARGE  Large sample properties of the regression depth induced median
 \\[4ex]
}

\author{ {\sc Yijun Zuo}\\[2ex]
         {\small {\em  Department of Statistics and Probability,  Michigan State University} }\\
         {\small East Lansing, MI 48824, USA} \\
         {\small zuo@msu.edu}\\[6ex]
     }
 \date{\today}
\maketitle

\vskip 3mm
{\small

\begin{abstract}
Notions of depth in regression have been introduced and studied in the literature.
 Regression depth (RD) of Rousseeuw and Hubert (1999), the most famous one, is a direct extension of Tukey location depth (Tukey, 1975) to regression.
   \vskip 3mm
Like its location counterpart, the most remarkable advantage of the notion of depth in regression is to directly introduce the maximum (or deepest) regression depth estimator (aka depth induced median) for regression parameters in a multi-dimensional setting.
Classical questions for the regression depth induced median include (i) is it a consistent estimator (or rather under what sufficient conditions, it is consistent)? and (ii) is there any limiting distribution?  \vs

Bai and He (1999) (BH99) pioneered an attempt to answer these questions. Under some stringent conditions on (i) the design points, (ii) the conditional distributions of $y$ given $\bs{x}_i$, and (iii) the error distributions, BH99 proved the strong consistency of the depth induced median. Under another set of conditions, BH99 showed the existence of the limiting distribution of the estimator.\vs

This article establishes the strong consistency of the depth induced median without any of the stringent conditions in BH99, and proves the existence of the limiting distribution of the estimator by sufficient conditions and an approach different from BH99.
\vskip 3mm
\vskip 3mm
\bigskip
\noindent{\bf AMS 2000 Classification:} Primary 62G09; Secondary
62G05, 62G15 62G20.
\bigskip
\par

\noindent{\bf Key words and phrase:} regression depth, maximum depth estimator, depth induced median, consistency, limiting distribution.
\bigskip
\par
\noindent {\bf Running title:} Asymptotic theorems for the regression depth median.
\end{abstract}
}

\section{Introduction}
Depth notions in location have received much attention in the literature. In fact, data depth and its applications remain one of the most active research topics in statistics in the last three decades.
Most favored notions of depth in location include (i) halfspace depth (HD)(Tukey, 1975, popularized by Donoho and Gasko, 1992),
(ii) simplicial depth (Liu, 1990),
and (iii) projection depth (PD) (Liu, 1992 and Zuo and Serfling, 2000, promoted by Zuo, 2003), among others.
\vskip 3mm

Depth notions in regression have also been sporadically proposed. Regression depth (RD) (Rousseeuw and Hubert, 1999) (RH99), 
the most famous one, is a direct extension of Tukey HD to regression.
Others include Carrizosa depth  (Carrizosa, 1996) and the projection regression depth induced from Marrona and Yohai, 1993 (MY93) and proposed in Zuo, 2018 (Z18). The latter turns out to be the extension of PD to regression.
\vs

One of the prominent advantages of depth notions
is that they can be directly employed to introduce median-type deepest estimating functionals (or estimators in the empirical case) for the location or regression parameters in a multi-dimensional setting based on a general min-max stratagem. The maximum (or deepest) regression depth estimator  serves as a \emph{robust} alternative to the classical least squares or least absolute deviations estimator of the unknown parameters in a general linear regression model. The latter can be expressed as \vspace*{-2mm}
\begin{eqnarray}
y&=&\mathbf{x}'\boldsymbol{\beta}+{{e}}, \label{eqn.model}
\end{eqnarray}
where  $'$ denotes the transpose of a vector, and random vector $\mathbf{x}=(x_1,\cdots, x_p)'$ and  parameter vector $\boldsymbol{\beta}$ are in $\R^p$ ($p\geq2$) and random variables $y$ and ${e}$ are in $\R^1$.
 If $\bs{\beta}=(\beta_0, \bs{\beta}'_1)'$ and $x_1=1$, then one has $y=\beta_0+\mb{x}'_1\bs{\beta}_1+{e}$, where $\mb{x}'_1=(x_2,\cdots, x_p) \in \R^{p-1}$.
Let $\mb{w}'=(1,\mb{x}'_1)$. Then $y=\mb{w}'\bs{\beta}+{e}$. We use this model or (\ref{eqn.model}) interchangeably depending on the context.
\vs
The maximum regression depth estimator possesses the outstanding robustness feature similar to the univariate location counterpart. 
 Indeed,
the maximum depth estimator induced from $RD$, could, asymptotically, resist up to $33\%$  contamination without breakdown,
in contrast to $0\%$ for the classical estimators 
(see Van Aelst and Rousseeuw, 2000) (VAR00). \vs

The asymptotics of the maximum regression depth estimator 
(denoted by $T^*_{RD}$, or $\bs{\beta}^*$)
 have been considered and established in Bai and He, 1999 (BH99). Under some stringent conditions on (i) the design points, (ii) the conditional distributions of $y$ given $\bs{x}_i$, and (iii) the error distributions, BH99 proved the strong consistency of the maximum depth estimator. Under another set of conditions, BH99 showed the existence of the limiting distribution of the estimator.
This article establishes the strong consistency of the maximum depth estimator without any of the stringent conditions in BH99 and proves the existence of the limiting distribution of the estimator with conditions and an approach very different from BH99.\vs
The rest of article is organized as follows: Section 2 introduces the RD of RH99 and presents pioneering
examples of the computation of the RD for population distributions. Section 3 summarizes the important results (from Z18) on the RD which are used in later sections. Section 4 establishes
the  strong and root-$n$  consistency  of the $T^*_{RD}$.
Section 5 is devoted to the establishment of limiting distribution of the $T^*_{RD}$, where the main tool is the Argmax continuous mapping theorem. Assumptions for the theorem to hold are verified via empirical process theory and especially stochastic equicontinuity and VC-classes of functions. The limiting distribution is characterized through an Argmax operation over the infimum of a function involving a Gaussian process.


\section{Regression depth of Rousseeuw and Hubert (1999)}
\begin{figure}[h!]
\vspace*{-6mm}
\includegraphics[width=\textwidth]{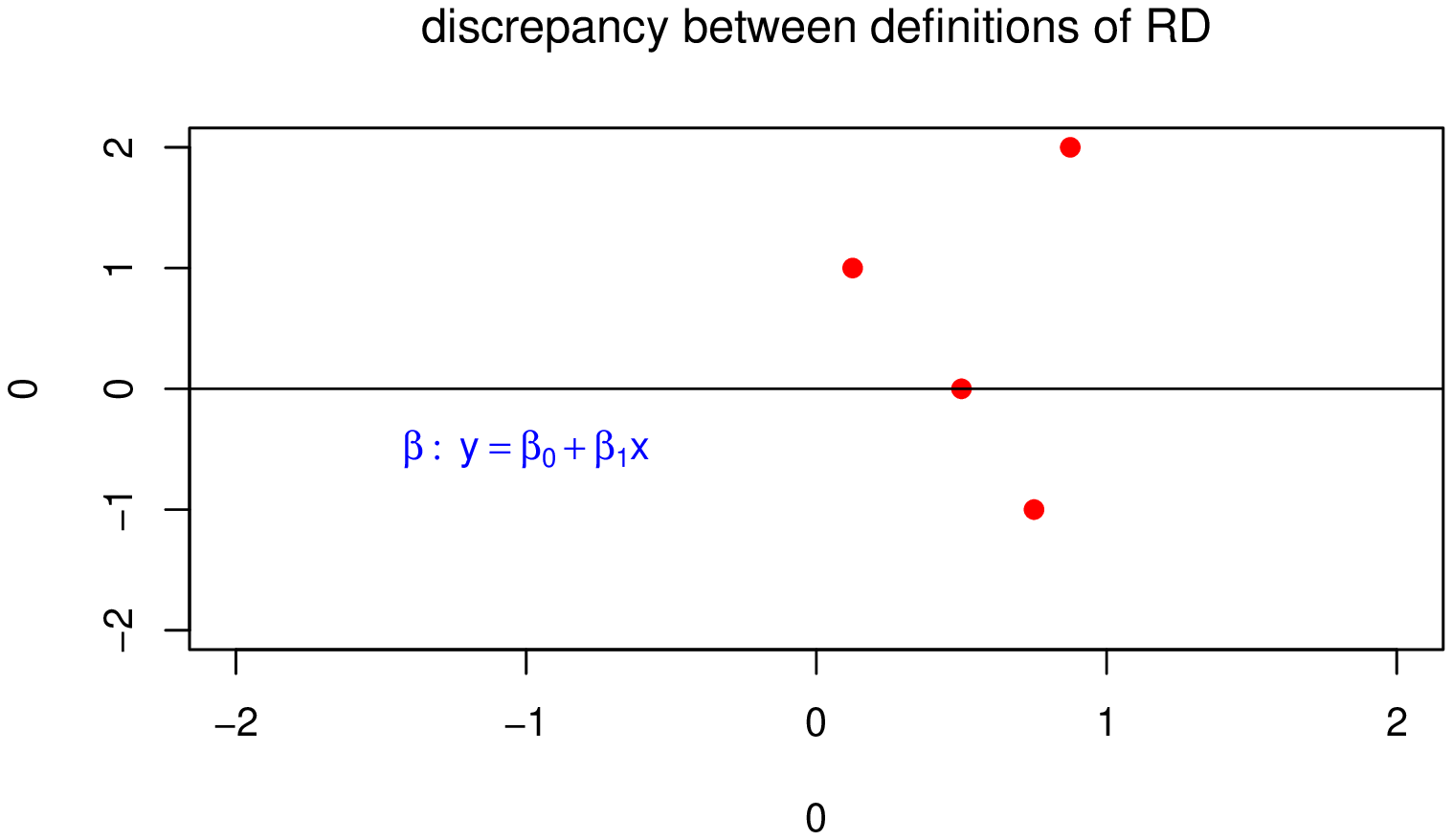}
\vspace*{-13mm}
 \caption{$\bs{\beta}=(0,0)'$, the horizontal candidate regression line. RH99 gives its RD=2 while RD$(\bs{\beta}, F_n)=1/2$ (Def. 2.1), RD of (\ref{bh99.eqn}) gives 1, whereas RD of (\ref{bh992.eqn}) yields $1.5$.
  }
 \label{fig:scatterplot-1}
\end{figure}

\vs
\tb{Definition 2.1} For any $\bs{\beta}$ and joint distribution $P$ of $(\mb{x}, y)$ in (\ref{eqn.model}), RH99 defined the regression depth of $\bs{\beta}$, denoted by RD$(\bs{\beta};P)$,
to be the minimum probability mass that needs to be passed when titling (the hyperplane induced from) $\bs{\beta}$ in any way until it is vertical. The maximum regression depth estimating functional $T^*_{RD}$ (also denoted by $\bs{\beta}^*$) is defined as
\be T^*_{RD}(P)=\arg\!\max_{\bs{\beta}\in\R^p}\mbox{RD}(\bs{\beta};P). \label{T-RD.eqn}
\ee
\noin
If there are several $\bs{\beta}$ that attain the maximum depth value on the right hand side (RHS) of (\ref{T-RD.eqn}), then the average of all those $\bs{\beta}$ is taken.\vs
The $\mbox{RD}(\bs{\beta};P)$ definition above is rather abstract and not easy to comprehend.
Some characterizations,  or equivalent definitions of RD$(\bs{\beta};P)$ are summarized below.
\emph{In the empirical case, the RD defined originally in RH99, divided by $n$, is identical to the following}.
\vs
\noin
\tb{Lemma 2.1}. The following statements for RD are equivalent.
\begin{itemize}
\item[] (i) [Z18]
 \be \mbox{RD}(\bs{\beta};P)=\inf_{\bs{\alpha} \in S(\bs{\beta})}P\left(|r(\bs{\beta})|\leq |r(\bs{\alpha})|\right),\ee
 where $S(\bs{\beta}):=\{\bs{\alpha}\in\R^p: ~ H_{\bs{\alpha}} \mbox{~intersects with~} H_{\bs{\beta}}\}$ for a given $\bs{\beta}$, $H_{\bs{\gamma}}$ denotes the unique hyperplane determined by $y=\mb{w'}\bs{\gamma}$, and $r(\bs{\gamma}):=y-(1,\mb{x'})\bs{\gamma}:=y-\mb{w}'\bs{\gamma}$.
\item[] (ii) [Z18]  \be
{\mbox{RD}}(\boldsymbol{\beta}; P)= \!\!\!\!
\inf_{\|\mb{v}_2\|=1,~v_1\in \R}\!\!\!\!E\Big(\mb{I}\big(r(\bs{\beta})(v_1,\mb{v}_2')\mb{w}\geq0\big)\Big)=\!\!
\inf_{\mb{v}\in\mbs^{p-1}}\!\!E\Big(\mb{I}\big(r(\bs{\beta})\mb{v}'\mb{w}\geq0\big)\Big), \label{eqn.rd}
\ee
where, $\mbs^{p-1}:=\{\mb{u}\in\R^p:\|\mb{u}\|=1\}$ and $\mb{I}(A)$ (and throughout) stands for the indicator function of the set $A$.
 \hfill \pend
\end{itemize}
Other characterizations are also given in the literature, e.g.,  in VAR00, in Rousseeuw and Struyf, 2004 (RS04),
 in Adrover, Maronna, and Yohai (2002), in Mizera (2002) (pages 1689-1690) and in BH99. The latter is specifically defined by
\be \mbox{RD}(\boldsymbol{\beta};{P_n})\!\!=\!\!
\inf_{\|\mb{u}\|=1,~ v \in \R} \min \bigg\{\sum_{i=1}^n\mb{I}(r_i(\boldsymbol{\beta})(\mathbf{u}'\mathbf{x}_i-v)>0), 
\sum_{i=1}^n\mb{I}(r_i(\boldsymbol{\beta})(\mathbf{u}'\mathbf{x}_i-v)<0) \bigg\}, \label{bh99.eqn}
\ee
\vs
Furthermore, BH99 depended solely on the following alternative definition:
\be \mbox{RD}(\boldsymbol{\beta};{P_n})=\frac{n}{2}+\frac{1}{2}\inf_{\bs{\gamma} \in\mbs^{p-1}}\sum_{i=1}^n\text{sgn}(y_i-\mb{w}_i\bs{\beta})\text{sgn}(\mb{w}'_i\bs{\gamma}).\label{bh992.eqn}
\ee
\vs
\noindent
\tb{Remarks 2.1}\vs
\noin
(I) If one assumes that $P(\mb{x'}\mb{u}=v)=0$ and $P(r(\bs{\beta})=0)=0$ for any $\mb{u}$, $v$, and $\bs{\beta}$, then Definition (\ref{bh99.eqn}) of BH99 above is \emph{identical} (a.s.) to the original definition of RH99. 
\vs
\noin
(II) Generally, definition (\ref{bh992.eqn})  is neither identical to the RD of RH99, nor to (\ref{bh99.eqn}).
For example,
 assume that we have four
sample points in $\R^2$ which could be regarded as from a continuous or discrete $\mb{Z}$ , $\mb{Z}_1=(\frac{1}{8},1)'; \mb{Z}_2=(\frac{4}{8},0)';\mb{Z}_3=(\frac{6}{8},-1)';\mb{Z}_4=(\frac{7}{8},2)'$. Then it is readily seen that for $\bs{\beta}=(0,0)'$, RH99 gives RD$=2$,
which  divided by $n=4$ leads to 
${1}/{2}$ (identical to Def. 2.1), (\ref{bh99.eqn}) gives 1 whereas (\ref{bh992.eqn}) yields $1.5$ (see Figure \ref{fig:scatterplot-1} for the scatterplot and the line).
\hfill \pend


\vs
For empirical distributions ($P=P_n$),  computing  RD$(\bs{\beta},P)$ is quite straightforward and examples have been given in RH99.
For a general distribution (probability measure) $P$,   concrete examples of expression of RD$(\bs{\beta},P)$ are not yet given in the literature herebefore.
For special classes of distributions, however, one could derive the explicit expression for RD$(\bs{\beta},P)$.
In the examples below, for simplicity, we again confine our attention to the case $p=2$. That is, we have a simple  linear regression model $y=\beta_0+\beta_1 x+e$.
\bec
\enc
\vspace*{-12mm}
 \vs
\noin
\tb{Example 2.1}
A random vector $\mb{X}\in\R^p$ is said to be \emph{elliptically distributed}, denoted by $\mb{X} \sim E(h; \bs{\mu}, \bs{\Sigma})$, if its density is
 of the form
\be
f(\mb{x}) = c|\bs{\Sigma}|^{-1/2}h\left((\mb{x} -\bs{\mu})'\bs{\Sigma}^{-1}(\mb{x} -\bs{\mu})\right),~ \mb{x}\in\R^p,  \label{ell.eqn}
\ee
where $c$ is the normalizing constant,  $\bs{\mu}$ is the coordinate-wise median vector (or the mean vector if it exists), $\bs{\Sigma}$ is a positive definite matrix which is proportional to the covariance matrix if it exists. 
Generally $h$ is a known function.
A straightforward transformation such as $\mb{Z}=\bs{\Sigma}^{-1/2}(\mb{X}-\bs{\mu})$ leads to $\mb{Z}\sim E(h; \bs{0},\bs{I_p})$.
\vs
To seek concrete expression for RD$(\bs{\beta},P)$ and for the simplicity 
 we restrict to the case $h(x)=\exp(-x^2/2)$, i.e., the bivariate normal class (laplace, logistic and t classes could be treated similarly). Namely, we have $(x, y)\sim N_2(\bs{\mu},\bs{\Sigma})$. After applying the transformation above, we can assume without loss of generality (w.l.o.g.) that
$(x, y)\sim N_2(\bs{0},I_2)$, where $I_2$ is a 2 by 2 identity matrix. For any $\bs{\beta}=(\beta_0, \beta_1)'$,  by the invariance of RD (see Z18 and Section 3), we can consider the depth of $\bs{\beta}$ w.r.t. the $P$ that corresponds to the   $N_2(\bs{0},I_2)$.
\bec
\begin{figure}[h!]
\vspace*{-7mm}
\includegraphics[width=\textwidth]{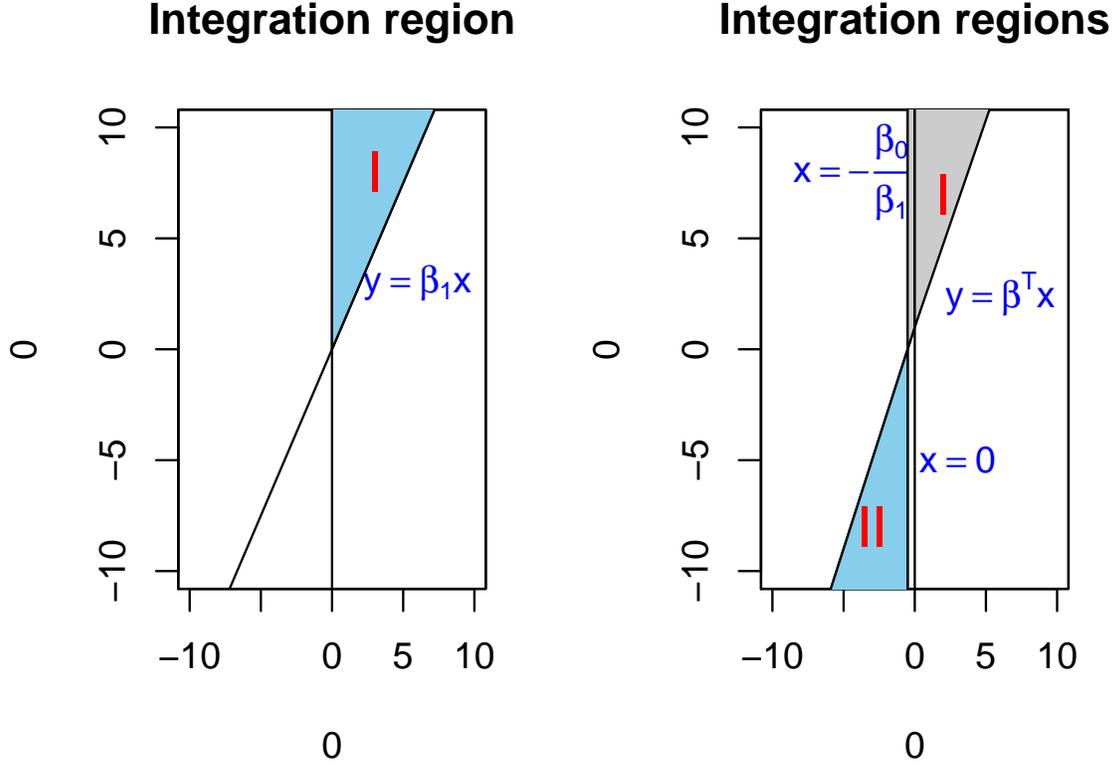}
\vspace*{-10mm}
 \caption{Integration regions. Left: for the region in (ii). Right: for the regions in (iii).
  }
 \label{fig:scatterplot}
\end{figure}
\enc
\vspace*{-13mm}

\begin{itemize}
\item[] (i) $\bs{\beta}=(0,0)'$, then the regression line is $y=0$, and $\text{RD}(\bs{\beta};P)=1/2$.

\item[] (ii) $\beta_0=0$ and $\beta_1>0$ (the case $\beta_1< 0$ can be discussed similarly). Denote the region bounded by
the regression line $y=\beta_1 x$ and the positive $y$-axis as I (see the left side of Fig. \ref{fig:scatterplot}), 
then  it is readily seen that
\be
\text{RD}(\bs{\beta};N(\mb{0};I_2))=2{P((y,x)\in \text{I})}\\[1ex] 
=1-2\int_0^{\infty}\Phi(\beta_1 x)d\Phi(x),
\ee
where 
 $\Phi(x)$ is the standard normal cumulative distribution function.
\vs
\item[] (iii) $\beta_0>0$ and $\beta_1>0$ (the case $\beta_0>0$ and $\beta_1 < 0$ and the cases where $\beta_0<0$   can be treated similarly).
Denote the region formed by the line with positive $y$ part of the vertical line $x=-\beta_0/\beta_1$,$\{x\geq -\beta_0/\beta_1, y\geq \beta_0+\beta_1 x\}$ as I and with negative $y$ part of the vertical line $x=-\beta_0/\beta_1$ $\{x\leq -\beta_0/\beta_1,y\leq \beta_0+\beta_1 x\}$ as II (see the right side of Fig.\ref{fig:scatterplot}), then it is readily seen
\bee
\text{RD}(\bs{\beta};N(\mb{0};I_2))={P((y,x)\in \text{I})+P((y,x)\in \text{II})}&&\\[1ex]
=1-2\Phi(-\beta_0/\beta_1)+\int_{-\infty}^{-\beta_0/\beta_1}\Phi\left(\beta_0+\beta_1 x\right)d\Phi(x) 
& & -\int^{\infty}_{-\beta_0/\beta_1}\Phi\left(\beta_0+\beta_1 x\right)d\Phi(x).
\ene
\item[] (iv) $\beta_0>0$ and $\beta_1=0$ (the case $\beta_0<0$ and $\beta_1=0$ can be handled similarly).
Denote the region  
formed by the line $y=\beta_0$ and the part of the positive $y$-axis $\{y\geq \beta_0\}$ as I
then it is readily seen that \vspace*{-2mm}
\bee
\text{RD}(\bs{\beta};N(\mb{0};I_2))&=& P((y,x)\in \text{I})=1/2-\Phi(\beta_0).
\ene
\hfill \pend
\end{itemize}
\noin
\tb{Example 2.2} Assume that $(x,y)$ is uniformly distributed over a unit circle centered at (0,0). By invariance of depth, this will cover a class of distributions of $A(x,y)+b$ for any nonsingular $A \in \R^{2\times2}$ and $b\in\R^2$.
\begin{itemize}
\item[] (i) $\bs{\beta}=(0,0)'$, then the regression line is  $y=0$, and $\text{RD}(\bs{\beta};P)=1/2$.
\item[] (ii) $\beta_0=0$ and $\beta_1>0$ ($\beta_1< 0$ can be treated similarly). Denote the region bounded by the regression line $y=\beta_1 x$ and the positive $y$ axis as I, then it is readily seen that
    \be
    \text{RD}(\bs{\beta};P)= 2P((x,y)\in \text{I}) = 1/2-\frac{|\arctan(\beta_1)|}{\pi}.
    \ee

\item[] (iii) $\beta_0>0$ and $\beta_1\geq  0$ (the cases where $(\beta_0>0, \beta_1 < 0)$  or $ (\beta_0<0$, $\beta_1\geq (\text{or} <)~0)$ can be dealt with similarly) and $\Delta=1+\beta_1^2-\beta_0^2> 0$. That is, the regression line intercepts the unit circle at two points $x_{\pm}$, where $x_{\pm}=\frac{-\beta_0\beta_1\pm \sqrt{1+\beta_1^2-\beta_0^2}}{1+\beta_1^2}$. \vs

    ~~~(a) Assume that both intersection points have positive $y$ coordinate.
     Denote the region formed by the regression line and the circle between the vertical lines $x=x_{-}$ and $x=x_{+}$, ($\beta_0+\beta_1 x\leq y\leq \sqrt{1-x^2}$) as I. Then it is readily seen that
      \bee
      \text{RD}(\bs{\beta};P)=P((x,y)\in \text{I})
      =\int^{x_{+}}_{x_{-}}\left(\sqrt{1-x^2}-(\beta_0+\beta_1x)\right)dx
      = g(x_{+})-g(x_{-})&&,
      \ene
 where $g(x)=\frac{1}{2}\left(x\sqrt{1-x^2}+\arcsin(x)\right)-\left(\beta_0x+\frac{1}{2}\beta_1x^2\right):=g_1(x)-g_2(x)$, 

      \vs
      ~~~(b) Assume that the $y$ coordinates of the two intersection points have different signs. The latter implies that $\beta_1\neq 0$. Denote  the region formed by the regression line and the circle  and the positive (negative) $y$-part of vertical line $x=-\beta_0/\beta_1$ as I (II).
      Then it is readily seen that
   \bee
    \text{RD}(\bs{\beta};P)=P((x,y)\in \text{I})+P((x,y)\in \text{II})
    = g(x_{+})-g(x_{-})+2g_2(-\beta_0/\beta_1)-2g_2(x_{-}). & &
   \ene


 \item[] (iv) In all other cases,   $\text{RD}(\bs{\beta};P)=0$. \hfill \pend
\end{itemize}
\vs
\noin
\tb{Remarks 2.1}\vs
(I) From the examples, it is readily seen that maximum value of RD is $1/2$ (in fact, $1/2$ is the maximum possible depth value in many cases, see RH99). Furthermore, the point $\bs{\beta}=(0,0)'$ is the unique point that attains the maximum depth value in both examples.
\vs
(II) According to RS04, we say $F_{\mb{Z}}$ is \emph{regression symmetric }about $\bs{\beta}^*=(0,0)'$ in these examples, where $\mb{Z}:=(x, y)$. We also have a unique T$^*_{RD}$ or $\bs{\beta}^*$ in both cases.
\hfill \pend

\section{Preliminary results}
\vs
\noin
A  regression depth functional $D$ is said to be regression, scale and affine \emph{invariant} w.r.t. a given $F_{(\mathbf{x}, y)}$ if and only if (iff),  respectively,
 $\displaystyle
 D(\boldsymbol{\beta}+\mb{b};~F_{(\mathbf{x}, ~y+\mathbf{x}'\mb{b})})=D(\boldsymbol{\beta};~F_{(\mathbf{x}, ~y)}), ~\forall ~\mb{b} \in \R^{p}$;~~ $ \displaystyle \!\!D(s\boldsymbol{\beta};~F_{(\mathbf{x},~sy)})=D(\boldsymbol{\beta};~F_{(\mathbf{x}, ~y)}),  ~\forall s (\neq 0)\in \R;~~$ $\displaystyle
\!\! D(A^{-1}\boldsymbol{\beta};F_{(A'\mathbf{x}, ~y)})=D(\boldsymbol{\beta}; F_{(\mathbf{x},~ y)}),\forall A_{p\times p},$
$\mbox{a nonsingular matrix}$. \vs
A regression estimating functional $T(\cdot)$ is said to be {regression}, {scale}, and {affine} \emph{equivariant} iff, respectively,
 $\displaystyle
 T(F_{(\mathbf{x}, ~y+\mathbf{x}'b)})=T(F_{(\mathbf{x}, ~y)})+b, ~\forall ~b \in \R^p;$  $\displaystyle
 T(F_{(\mathbf{x}, ~sy)})=s T(F_{(\mathbf{x}, ~y)}), ~\forall ~s\in \R;$  $\displaystyle
 T(F_{(A'\mathbf{x}, ~y)})=A^{-1}T(F_{(\mathbf{x}, ~y)}),~ ~ \forall \mbox{~ nonsingular  $A\in \R^{p\times p}$}.
$
\vs
We now summarize some preliminary results on the RD and its induced maximum depth estimating functional.
$F_{\mb{Z}}$ and P are used interchangeably and $\mb{Z}:=(\mb{x}, ~y)$. 
\vs
\noin
\tb{Lemma 3.1} [Zuo (2018)] \vs
\begin{itemize}
\item[] (i) $\text{RD}(\bs{\beta}; F_{\mb{Z}})$ is regression, scale and affine invariant  and hence $T_{RD}^*(F_{\mb{Z}})$
is regression, scale and affine equivariant. Furthermore, $\text{RD}(\bs{\beta}; F_{\mb{Z}})\to 0$ as $\|\bs{\beta}\| \to \infty$,
if $P(H_v)=0$ for any vertical hyperplane $H_v$.
\item[] (ii) $\text{RD}(\bs{\beta}; P)$ is upper-semicontinuous and concave (in $\bs{\beta}$), and continuous in $\bs{\beta}$ if $P$ has a density.

\item[] (iii) $\sup_{\bs{\beta}\in\R^p}|\text{RD}(\bs{\beta};F^n_{\mb{Z}})-\text{RD}(\bs{\beta};F_{\mb{Z}})|\to 0$ almost surely (a.s.) as $n\to \infty$, where $F^n_{\mb{Z}}$ is the empirical version of the distribution $F_{\mb{Z}}$. 
    \hfill \pend
\end{itemize}
\vs
In the sequel, we assume that there exists a unique point  $T_{RD}^*(F_{\mb{Z}})$ (or $\bs{\beta^*}$, a generic notation for the maximum regression depth point, whereas $\bs{\eta}$ is also used for the maximum (location or regression) depth point later)  that maximizes the underlying regression depth RD. In virtue of (i) of Lemma 3.1, 
one can assume, w.l.o.g., that $T_{RD}^*(F_{\mb{Z}})=\mb{0}$. \vs
Uniqueness is guaranteed if $F_{\mb{Z}}$ has a strictly positive density and is regression symmetric about a point $\bs{\beta}$ ($F_{\mb{Z}}$ is \emph{regression symmetric} about $\bs{\theta}$ if $P(\mb{x}\in B, r(\bs{\theta})>0)= P(\mb{x}\in B, r(\bs{\theta})<0)$ for any Borel set $B \in\R^{p-1}$, see RS04).

\vs
\section{Consistency}
For a general regression depth functional  $D(\bs{\beta};F_{\mb{Z}})$,
let $\bs{\beta^*}(F_{\mb{Z}})=\arg\max_{\bs{\beta}\in\R^p}D(\bs{\beta};F_{\mb{Z}})$, then
$\bs{\beta^*_n}:=\bs{\beta^*}(F^n_{\mb{Z}})$ is a natural estimator of $\bs{\beta^*}$, the maximum regression depth functional. 
\vs
Is  $\bs{\beta^*_n}$ a consistent estimator? This is a very typical question asked in statistics and the argument (or answer) for it is also very standard, almost to the point of clich\'{e} as Kim and Pollard (1990) (KP90) have commented.
\vs
Let us first deal with the problem in a more general setting. Let $M_n$  be stochastic processes indexed by a metric space $\Theta$ of $\bs{\theta}$, and $M\!\!:$ $\Theta \to\R$ be a deterministic function of $\bs{\theta}$ which attains its maximum at a point $\bs{\theta}_0$.\vs
 The sufficient conditions  for the consistency of this type of problem were given in Van Der Vaart (1998) (VDV98) and Van Der Vaart and Wellner (1996) (VW96) and are listed below:\vs
\begin{itemize}
\item[]\tb{C1:} $\sup_{\bs{\theta}\in\Theta}|M_n(\bs{\theta})-M(\bs{\theta})|=o_p(1)$;
\vskip 3mm
\item[]\tb{C2:} $\sup_{~\{\bs{\theta}:~d(\bs{\theta},\bs{\theta_0})\geq \delta\}} M(\bs{\theta})< M(\bs{\theta_0})$, for any $\delta>0$ and the metric $d$ on $\Theta$;
\vskip 3mm
\hspace*{-10mm} Then any sequence $\bs{\theta}_n$ is consistent for $\bs{\theta}_0$ providing  that it satisfies
\item[]\tb{C3:} $M_n(\bs{\theta_n})\geq M_n(\bs{\theta_0})-o_p(1)$.
\vskip 3mm

\end{itemize}

\vs
\noin
\tb{Lemma 4.1} [Th. 5.7, VDV98] If \tb{C1} and \tb{C2} hold, then any $\bs{\theta}_n$   satisfying \tb{C3} is consistent for $\bs{\theta}_0$.
\vs
\vs
\noin
\tb{Remarks 4.1}
\begin{itemize}
\vs
\item[] (I) \tb{C1} basically requires that the $M_n(\bs{\theta})$  converges   to $M(\bs{\theta})$ in probability  uniformly in $\bs{\theta}$. For the depth process RD$(\bs{\beta}; F^n_{\mb{Z}})$ and RD$(\bs{\beta}; F_{\mb{Z}})$, it holds true 
    (the convergence here is actually almost surely (a.s.)
     and uniformly in $\bs{\beta}$).
\vs
\item[] (II) \tb{C2} essentially demands that the unique maximizer $\bs{\theta_0}$ is well separated. This holds true 
    as long as $D(\bs{\beta};F_{\mb{Z}})$ is upper semi-continuous and vanishing at infinity, and $\bs{\theta_0}$ is unique (see, Lemma 4.2 below). Therefore, it holds for RD in light of Lemma 3.1.
\vs
\item[] (III) \tb{C3} asks that $\bs{\theta_n}$
 is very close to $\bs{\theta}_0$ in the sense that the difference of images of the two at $M_n$ is within $o_p(1)$.
In KP90 and VW96 a stronger version of \tb{C3} is required:
\vs
$\displaystyle \mbox{\tb{C3}'}: ~~ M_n(\bs{\theta_n})\geq \sup_{\theta\in\Theta} M_n(\bs{\theta})-o_p(1),$\vs
 which implies \tb{C3}. This strong version mandates that $\bs{\theta_n}$  \emph{nearly} maximizes $M_n(\bs{\theta})$.
The maximum regression depth estimator $\bs{\beta^*_n}(:=\bs{\theta_n})$ is defined to be the maximizer of $M_n(\bs{\theta}):=D(\bs{\beta};F^n_{\mb{Z}})$, hence \tb{C3}' (and thus \tb{C3}) holds automatically. \hfill \pend

\end{itemize}
\noin
In light of above, 
$\bs{\beta^*_n}$ induced from RD is consistent for $\bs{\beta^*}$. But, we have more.
\vs
\noindent
\tb{Theorem 4.1} The maximum regression depth estimator $\bs{\beta^*_n}$ induced from RD 
 is \emph{strongly} consistent for $\bs{\beta^*}$ (i.e., $\bs{\beta^*_n}-\bs{\beta^*}=o(1)$ a.s.) provided that $\bs{\beta}^*$ is unique.
\vs
\noindent
\tb{Proof:} 
The proof for the consistency of Lemma 4.1 could be easily extended to the strong consistency with a strengthened version  of \tb{C1}
 \vs
$\displaystyle ~~~~\mbox{\tb{C1}':}~~\sup_{\bs{\theta}\in\Theta}|M_n(\bs{\theta})-M(\bs{\theta})|=o(1)$ a.s.\vs
In the light of the proof of Lemma 4.1, we need only  verify the sufficient conditions \tb{C1}' and \tb{C2-C3}.  By (III) of Remark 4.1,
\tb{C3} holds automatically, so  we need to verify \tb{C1}' and \tb{C2}.
\tb{C1}' has been given in Lemma 3.1 for RD. So the only item left is to verify  \tb{C2} for
RD which is guaranteed by Lemma 4.2 below. 
\hfill \pend
\vs
\noindent
\tb{Lemma 4.2} Assume that a general (location or regression) depth $D(\bs{\beta}; F_{\mb{Z}})$ is upper semi-continuous in $\bs{\beta}$ and vanishing when $\|\bs{\beta}\|\to \infty$.
 Let $\bs{\eta}\in\R^p$ be the unique point with $\bs{\eta}=\arg\max_{\bs{\beta}\in\R^p}D(\bs{\beta}; F_{\mb{Z}})$ and $D(\bs{\eta}; F_{\mb{Z}})>0$.
Then for any $\ep>0$, $\sup_{\bs{\beta}\in N^c_{\ep}(\bs{\eta})}D(\bs{\beta}; F_{\mb{Z}}) < D(\bs{\eta}; F_{\mb{Z}})$, where $N^c_{\ep}(\bs{\eta})=\{\bs{\beta}\in\R^P: \|\bs{\beta}-\bs{\eta}\|\geq \ep\}$ and ``c" stands for ``complement" of a set.\vskip 3mm
\noin
\tb{Proof}:  Assume conversely that $\sup_{\bs{\beta}\in N^c_{\ep}(\bs{\eta})}D(\bs{\beta}; F_{\mb{Z}})=D(\bs{\eta}; F_{\mb{Z}})$. Then by the given conditions, there is a sequence of \emph{bounded} $\bs{\beta_j}$ ($j=0,1,\cdots$) in $N^c_\ep(\bs{\eta})$ such that   $\bs{\beta_j}\to \bs{\beta_0}\in N^c_\ep(\bs{\eta})$ and $D(\bs{\beta_j}; F_{\mb{Z}})\to D(\bs{\eta}; F_{\mb{Z}})$ as $j\to \infty$.  Note that $D(\bs{\eta}; F_{\mb{Z}})>D(\bs{\beta_0}; F_{\mb{Z}})$. The upper-semicontinuity of $D(\cdot; F_{\mb{Z}})$ now leads to a contradiction: for sufficiently large $j$, $D(\bs{\beta_j}; F_{\mb{Z}})\leq (D(\bs{\eta}; F_{\mb{Z}}) +D(\bs{\beta_0};F_{\mb{Z}}))/2< D(\bs{\eta}; F_{\mb{Z}})$.  This completes the proof. \hfill \pend
\vs
\noindent
\tb{Remarks 4.2}
\vs
\bi
\item[] (I) For RD, the sufficient conditions in the Lemma are all satisfied in virtue of Lemma 3.1 and the uniqueness of  $\bs{\eta}=\bs{\beta}^*$ is guaranteed for special $F_{\mb{Z}}$ (see Section 3).
\item[] (II) Besides the necessary uniqueness assumption of $\bs{\beta}^*$ here, BH99 under additional stringent conditions  (see their D1-D4)
on (i) design points $\mb{x}_i$, (ii) the conditional distributions of $y$ given $\mb{x_i}$, and (iii) the distributions of error $e_i$, proved
the strong consistency with a very different and unnecessarily complicated approach.\vs (III) However, if i.i.d $x_i$ is from a univariate Cauchy distribution or any other heavy- tailed ones, then all D1-D3 do not hold,  Theorem 2.1 of BH99 is not applicable and no strong consistency result can be obtained meanwhile one can get the result via Theorem 4.1 above, nevertheless. This is the merit and necessity of Theorem 4.1.  \hfill \pend

\ei

\vs
\noindent
With the establishment of  strong  consistency, one naturally  wonders about the rate of convergence of the maximum regression depth estimator. Does it possess  root-$n$ consistency?
To answer the question,  we need a stronger version of \tb{C2} for a general depth notion $D$.
\vs
\bi
\item[] \tb{C2}': For each small enough positive $\delta$, there exists a positive constant $\kappa$ such that
 $\sup_{~\|\bs{\beta}-\bs{\eta}\|\geq \delta} D(\bs{\beta};P)< \alpha^*-\kappa \delta,$
  where $\alpha^*:=D(\bs{\eta};P)=\sup_{\bs{\beta}\in\R^p}D(\bs{\beta};P)$.
\ei
\noin
\tb{Remarks 4.3:}\vs
(I)  When $D$ in \tb{C2}' is RD, Lemma 4.2 provides a choice for the individual $\kappa$ for every $\delta$. But \tb{C2}' requires more. 
In the following we provide sufficient conditions for \tb{C2}' to hold. 
\vs
(II) (i) $P$ has a density; (ii) $h(\bs{\beta}, \mb{v}):=E\left(\mathbf{I}(r(\bs{\beta}) \mb{v}'\mb{w}\geq 0)\right)$ is differentiable in $\bs{\beta}\in N_{\bs{\eta}}$ for a given $\mb{v} \in\mbs^{p-1}$, where $N_{\bs{\eta}}$ is a small neighborhood of $\bs{\eta}$;
      (iii) the directional derivative of $h$ along  $\mb{u}\in \mbs^{p-1}$ at $\bs{\beta}\in N_{\bs{\eta}}$
: $D_{\mb{u}}h(\bs{\beta},\mb{v}):=\nabla h(\bs{\beta}, \mb{v})\cdot \mb{u}$ is continuous in $\mb{u}$ and positive uniformly in $\mb{u}$, where $\nabla$ is the vector differential operator and ``$\cdot$" stands for the inner product; (iv) the Hessian matrix of $h$ has a positive eigenvalues uniformly for $\bs{\beta}$ over $N_{\bs{\eta}}$.
\hfill \pend
\vs
Let $D(\bs{\beta};P)$ be any \emph{regression} (or even \emph{location}) depth functional for $\bs{\beta}\in\R^p$. We have the following general result for
$\bs{\beta^*_n}=\arg\max_{\bs{\beta}\in\R^p}D(\bs{\beta}; P_n)$ and $\bs{\beta^*}=\arg\max_{\bs{\beta}\in\R^p}D(\bs{\beta}; P)$:
\vs
\noindent
\tb{Lemma 4.3} Let $D(\bs{\beta};P)$ be a general depth notion. 
 If
(i) $\sup_{\bs{\beta}\in\R^p}|D(\bs{\beta};P_n)-D(\bs{\beta};P)|=O_p(n^{-1/2})$ and (ii) \tb{C2}' holds, 
then $\bs{\beta}^*_n-\bs{\beta^*}=O_p(n^{-1/2})$.
\vs
\noin
\tb{Proof}:~~
Denote  $\Delta_n:= \sup_{\bs{\beta}\in\R^p}|D(\bs{\beta};P_n)-D(\bs{\beta};P)|$.
 Let $\delta=2\Delta_n/\kappa$. In light of \tb{C2}', we have for every $n$
\bee
\sup_{\bs{\beta}:\|\bs{\beta}-\bs{\beta}^*\|\geq \delta}\!\!\!D(\bs{\beta};P_n)
\!\!\!&\leq&\!\!\! \sup_{\bs{\beta}:\|\bs{\beta}-\bs{\beta}^*\|\geq \delta}\!\!\!|D(\bs{\beta};P_n)-D(\bs{\beta};P)| +\!\!\!\sup_{\bs{\beta}:\|\bs{\beta}-\bs{\beta}^*\|\geq \delta}\!\!\!D(\bs{\beta};P)\\[1ex]
&<& \Delta_n+  \alpha^*-\kappa \delta 
=\alpha^*-\Delta_n,
\ene
which, in conjunction with the definition of $\bs{\beta^*_n}$, 
 implies that $\|\bs{\beta}^*_n-\bs{\beta}^*\|< \delta$.
\hfill \pend
\vs
\noin
\tb{Theorem 4.2} If \tb{(A0):}  $E|r(\bs{\beta})|^2$ and $E|\mb{v}'\mb{w}|^2$  exist 
uniformly in $\bs{\beta}\in \R^p$ and $\mb{v}\in\mbs^{p-1}$, then
(i) $\sup_{\bs{\beta} \in\R^p}|\mbox{RD}(\bs{\beta};P_n)-\mbox{RD}(\bs{\beta};P)|=O_p(n^{-\frac{1}{2}})$ and 
(ii)
   $\bs{\beta}_n^*-\bs{\beta}^*=O_p(n^{-\frac{1}{2}})$ 
   if \tb{C2}' holds. 
 \vs
 \noindent
\tb{Proof:}~~
  Write $f(y,\mb{w},\bs{\beta}, \mb{v}):=(y-
 \mb{w}'\bs{\beta})(\mb{v'}\mb{w})$, $\forall \bs{\beta}\in\R^p, \mb{v}\in\mbs^{p-1}$.  In  light of (ii) of Lemma 2.1,  RD$(\bs{\beta};P)=\inf_{\mb{v}\in\mbs^{p-1}}P\left(f(y,\mb{w}, \bs{\beta},\mb{v})\geq 0\right)$. Define a class of functions (for the notation convention, see p140 of Pollard,1984 (P84))
$
\mathscr{F}=\{(f(\cdot,\cdot,\bs{\beta},\mb{v}),~\forall~ \bs{\beta}\in\R^p,\text{~and~} \mb{v}\in\mbs^{p-1} \}.
$
Note that (see Z18)
\bee
&&\{f(y,\mb{w}, \bs{\beta}, \mb{v})\geq 0, \bs{\beta}\in \R^p, \mb{v}\in \mbs^{p-1} \}=
\Big\{\Big(\{(y-\bs{\beta}' \mb{w})\geq 0\}\cap\{\mb{u}' \mb{x}< v\}\Big)\\ &&\cup \Big(\{ (y-\bs{\beta}' \mb{w})<0\}\cap\{\mb{u}' \mb{x}\geq v\}\Big),   \mb{u}\in \mbs^{p-2}, v\in \R^1, \bs{\beta}\in \R^p\Big\}.
\ene
The RHS above is built up from sets of the form $\{g\geq 0\}$ with $g$ in the finite-dimensional
vector space of functions. By Lemmas II.28 and II.15 of P84, the class of graphs of functions in  $\mathscr{F}$
has polynomial discrimination 
 and  $\mathscr{F}$
  has VC subgraphs with a square integrable envelope (see II.5 of P84 and 2.6 of VW96 for discussions). By Corollary 3.2 of KP90,
we have that
\[
\sup_{\bs{\beta}\in\R^p,~ \mb{v}\in\mbs^{p-1}}\big|P_n(f(y,\mb{w},\bs{\beta},\mb{v})\geq 0)-P(f(y,\mb{w},\bs{\beta},\mb{v})\geq 0)\big|=O_p(n^{-1/2}).
\]
Thus  we have
 \bee
 \sup_{\bs{\beta}\in\R^p}
 \bigg| \inf_{\mb{v}\in\mbs^{p-1}}P_n(f(y,\mb{w},\bs{\beta},\mb{v})\geq 0)-\inf_{\mb{v}\in\mbs^{p-1}}P(f(y,\mb{w},\bs{\beta},\mb{v})\geq 0)\bigg|
 &\leq&\\[1ex]
\sup_{\bs{\beta}\in\R^p,~ \mb{v}\in\mbs^{p-1}}\big|P_n(f(y,\mb{w},\bs{\beta},\mb{v})\geq 0)-P(f(y,\mb{w},\bs{\beta},\mb{v})\geq 0)\big|
&=&O_p(n^{-1/2}),\\
\ene
where the inequality follows from the fact that $|\inf_A f-\inf_A g|\leq \sup_A|f-g|$. It follows that
$\sup_{\bs{\beta}\in\R^p}|\mbox{RD}(\bs{\beta};P_n)-\mbox{RD}(\bs{\beta};P)|=O_p(n^{-1/2})$, the first part of the theorem is obtained.
\vs
\vs
This first part, in conjunction with \tb{C2}' and Lemma 4.3, yields the desired second part of the theorem. That is,
$\bs{\beta}^*_n-\bs{\beta}^*=O_p(n^{-1/2})$.
\hfill \pend
\vs
\vs
\noindent
\tb{Remarks 4.4:}
\bi

\item[] (I) The approach of the first part of the proof could be extended for any depth notions that are defined based on sets that form a VC class such as the location counterpart, Tukey halfspace depth (HD), where one has a class of halfspaces, a VC class of sets.
    \vs
     That is, utilizing the approach, one can prove that the maximum Tukey location depth estimator (aka Tukey median) is root-$n$ consistent (uniformly tight) if \tb{C2}' holds for the HD. 
     For  the latter, a sufficient condition was given in Nolan (1999) ((ii) of Lemma 2), BH99 ((N2) in Theorem 4.1), and Mass\'{e} (2002) ((b) of Proposition 3.2 and Theorem 3.5). That is, the approach here covers the uniform tightness result in those papers.

     \vs
\item[] (II)  BH99 obtained the root-$n$ consistency for  $\bs{\beta^*_n}$ with a very different approach under more assumptions, such as their (D1)-(D4) and (C1), (C2), and (C3), on the random vector $\mb{x}$,  on the conditional distribution of $y$ given $\mb{x}$, and on the error distributions.
 \hfill \pend
 \ei
 \vs
\section{Limiting distribution}
With the root-$n$ consistency of the maximum regression depth estimator established, we are now in a position to address the natural question: does it have a limiting distribution?
\vs
Since the tool employed for establishing the limiting distributions is the \emph{Argmax} theorem, we first cite it below from VW96 (Theorem of 2.7 of KP90 is an earlier version).\vs
\noin
\tb{Lemma 5.1} [Th. 3.2.2, VW96, Argmax continuous mapping]  Let $M_n$, $M$ be stochastic processes
indexed by a metric space $S$ such that $M_n\xrightarrow{d} M$ in ${l}^{\infty}(K)$ for every compact
$K\subset S$. Suppose that almost all sample paths $\mb{s} \mapsto M(\mb{s})$ are upper semicontinuous
and possess a unique maximum at a random point $\widehat{\mb{s}}$, which, as a
random map into $S$, is tight. If the sequence $\widehat{\mb{s}}_n$ is uniformly tight and
satisfies $M_n(\widehat{\mb{s}}_n) \geq \sup_{\mb{s}} M_n(\mb{s})-o_p(1)$, then $\widehat{\mb{s}}_n\xrightarrow{d} \widehat{\mb{s}}$, where $\xrightarrow{d}$ stands for convergence in distribution. 
\hfill \pend
\vs
To establish the limiting distribution for $\widehat{\mb{s}}_n:=\sqrt{n}\bs{\beta^*_n}$, we need (A) to identify the processes $M_n$ and $M$ and show that $M_n\xrightarrow{d} M$  in $l^{\infty}(K)$ for any compact $K\subset \R^p$. (B) to show that almost all sample paths of $M(\mb{s})$ are upper semicontinuous and  possess a unique maximum at a random point $\widehat{\mb{s}}$, which is tight, and (C) to show that $\widehat{\mb{s}}_n$ is uniformly tight and $M_n(\widehat{\mb{s}}_n) \geq \sup_{s} M_n(\mb{s})-o_p(1)$.\vs

In virtue of Theorem 4.2, part of (C) already holds under certain conditions for $\widehat{\mb{s}}_n=\sqrt{n}\bs{\beta^*_n}$. So we need to verify the (A) and (B) and the second part of (C).
\vs
\vs
By (ii) of Lemma 2.1 and (\ref{T-RD.eqn}),
 we have that $$\bs{\beta^*}=\arg\!\max_{\bs{\beta}\in\R^p}\mbox{RD}(\bs{\beta};P)=
\arg\!\max_{\bs{\beta}\in\R^p}\inf_{\mb{v}\in\mbs^{p-1}}E\Big(\mathbf{I}\big(f(y,\mb{w},\bs{\beta},\mb{v})\geq0\big)\Big),$$
where $f(y,\mb{w},\bs{\beta},\mb{v})=(y-\mb{w'}\bs{\beta})\mb{v'}\mb{w}$ given in the proof of Theorem 4.2.
For a given $\bs{\beta}$ define
$$
{V}(\bs{\beta})=\{\mb{v}\in\mbs^{p-1}: \mbox{RD}(\bs{\beta};P)=P(f(y,\mb{w},\bs{\beta},\mb{v})\geq0)=\!\!\inf_{\mb{u}\in\mbs^{p-1}}\!\!P(f(y,\mb{w},\bs{\beta},\mb{u})\geq0)\},
$$
i.e.,
the collection of $\mb{v}$ at which $P(f (y,\mb{w},\bs{\beta},\mb{v})\geq 0)$ attains the infimum over $\mb{v}\in\mbs^{p-1}$.
\vs
Assume by Lemma 3.1 that $\bs{\beta^*}$ is $\mb{0}$.  Hereafter $\bs{\beta}$ is assumed to be in a small bounded neighborhood $\Theta$ of $\mb{0}$ by virtue of Theorem 4.1.
Assume for  $\mb{v}\in\mbs^{p-1}$ and  $\bs{\beta}\in \Theta$ that
\[ \hspace*{-3mm}
\tb{A1}:~~~~~P(f(y,\mb{w},\bs{\beta},\mb{v})\geq 0)=P(f(y,\mb{w},\mb{0},\mb{v})\geq 0)+\mb{g}(\mb{v})\cdot\bs{\beta}+o(\|\bs{\beta}\|), \hfill
\]
where $\mb{g}(\mb{v})$  is the $\nabla h(\bs{0}, \mb{v})$. The latter is defined in (ii) of (II) of Remarks 4.3. %
That is, the LHS permits a Taylor expansion at $\bs{\beta^*}=\mb{0}$.
Furthermore, 
\[
\hspace*{-53mm}\tb{A2}:~~~~~~~~~\hspace*{18mm} V(\mb{0})=\mbs^{p-1}. 
 \hspace{13mm}
\]
That is, along any direction $\mb{v}\in\mbs^{p-1}$, $P(f(y,\mb{w},\mb{0},\mb{v})\geq 0)=\alpha^*:=\mbox{RD}(\bs{\beta^*},P)$. 
And
\[
\hspace*{-5mm}\tb{A3}:~~~~~~~~~\mbox{$\mb{g}(\mb{v})$ is continuous in $\mb{v}$ and} \sup_{\mb{v}\in V(\mb{0})}\|\mb{g}(\mb{v})\|<\infty.
\hspace{20mm}
\]
That is, $\mb{g}(\mb{v})$  is  uniformly bounded  over $V(\mb{0})$.
\vs
\noindent
\tb{Theorem 5.1} If \tb{C2}' and \tb{A0}-\tb{A3} hold, 
then  for $\bs{\beta^*_n}$ induced from RD,  as $n\to \infty$, $$\sqrt{n}(\bs{\beta^*_n}-\bs{\beta^*})\xrightarrow{~~~d~~~}
\arg\!\max_{\mb{s}} \inf_{\mb{v}\in {V}(\mb{0})}\{{E_P}(f(y,\mb{w},\mb{0},\mb{v})\geq 0)+\mb{g}(\mb{v})\cdot\mb{s}\},$$
where ${E_P}$ is the limit of the empirical process $E_n=\sqrt{n}(P_n-P)$ in $l^{\infty}(\mathscr{F})$, a $P$-Brownian bridge (see Def. VII. 14 of P84), and $\mathscr{F}=\{f(\cdot,\cdot,\mb{0},\mb{v})\geq 0, \mb{v}\in V(\mb{0})\}$.
\vs
\vs
\noin
\tb{Proof}: ~~
\tb{C2}' guarantees the uniqueness of $\bs{\beta}^*(:=\bs{\eta})$, which can be assumed, w.l.o.g., to be $\mb{0}$.
 (ii) of Lemma 2.1 yields
$$\mbox{RD}(\bs{\beta}; P)=\inf_{\mb{v}\in\mbs^{p-1}}E(\mathbf{I}({(y-\mb{w}'\bs{\beta})\mb{v}'\mb{w}\geq 0}))=\inf_{\mb{v}\in\mbs^{p-1}}P(f(y,\mb{w},\bs{\beta},\mb{v})\geq 0).$$
Note that
\be n^{1/2}\bs{\beta}^*_n=n^{1/2}\arg\!\max_{\bs{\beta}\in\R^p}\inf_{\mb{v}\in\mbs^{p-1}}P_n\left(f(y,\mb{w},\bs{\beta},\mb{v})\geq 0\right).\label{betan*.eqn}
\ee
Hence for any compact $K$ and $\mb{s}\in K\subset \R^p$ and sufficiently large $n$, 
\begin{align}
&n^{1/2}P_n(f(y,\mb{w},\mb{s}/n^{1/2},\mb{v})\geq 0)\!=\!n^{1/2}P(f(y,\mb{w},\mb{s}/n^{1/2},\mb{v})\geq 0) 
+E_n(f(y,\mb{w},\mb{s}/n^{1/2},\mb{v})\geq 0) \nonumber\\[1ex]
                            & \!=\!n^{1/2}P(f(y,\mb{w},\mb{s}/n^{1/2},\mb{v})\geq 0)+E_n(f(y,\mb{w},\mb{0},\mb{v})\geq 0)+o_p(1)\nonumber\\[1ex]
                            &\! =\!n^{1/2} P(f(y,\mb{w},\mb{0},\mb{v})\geq 0)+ \mb{g}(\mb{v})\cdot\mb{s}+o(\|\mb{s}\|)
                              +E_n(f(y,\mb{w},\mb{0},\mb{v})\geq 0)+o_p(1), \label{taylor.eqn}
\end{align}
where 
 the second equality follows from the stochastic equicontinuity Lemma VII. 15 of P84 (see, Lemma II.18; Example II.26; Lemma II. 28; Lemmas II.25, II.36 and Example VII.18 of P84),
 the last equality follows from the \tb{A1}.
Then we can define that
\be M_n(\mb{s}):=n^{1/2}\inf_{\mb{v}\in\mbs^{p-1}}P_n(f(y,\mb{w},\mb{s}n^{-1/2},\mb{v}))-n^{1/2}\alpha^*,\label{mns.eqn}
\ee
where $\alpha^*=\mbox{RD}(\bs{\beta^*};P)$.
Note that by (\ref{betan*.eqn}), it is readily seen that $\widehat{\mb{s}}_n:=n^{1/2}\bs{\beta}^*_n$  maximizes $M_n(\mb{s})$ and is uniformly tight in virtue of Theorem 4.2,  therefore (C) is completely verified.
\vs
Now we need to verify (A) and (B) for
 \be  M(\mb{s}):=\inf_{\mb{v}\in {V}(\mb{0})}\{{E_P}(f(y,\mb{w},\mb{0},\mb{v})\geq 0)+\mb{g}(\mb{v})\cdot\mb{s}\}.\label{ms.eqn} \ee

\noin
We first establish some lemmas to fulfil the task above. \vs
\noin
\tb{Lemma 5.2} In light of \tb{A2} and \tb{A3},
 \bi
\item[] \tb{R1:} The sample path of $M(\mb{s})$ is continuous in $\mb{s}$ a.s., and furthermore $M(\mb{s})\to -\infty$ as $\|\mb{s}\|\to \infty$ a.s.;
~~\tb{R2}: $M(\mb{s})$ is concave in $\mb{s}$ a.s..
 \ei
 \noindent
 \tb{Proof}: Write $M(\mb{s},\mb{v})={E_P}(f(y,\mb{w},\mb{0},\mb{v})\geq 0)+\mb{g}(\mb{v}) \cdot \mb{s}$. The continuity and concavity of $M(\mb{s},\mb{v})$ in $\mb{s}$ is obvious. The assertion on $M(\mb{s})$ follows since the infimum preserves these properties. We need to show the second part of \tb{R1}.\vs

 First by the compactness of $V(\mb{0})$, the continuity and boundedness of $g(\mb{v})$ over $V(\mb{0})$, for an arbitrary $\mb{s}$, there is a $\mb{v}_0\in V(\mb{0})$ such that
 \be\inf_{\mb{v}\in V(\mb{0})}\mb{g}(\mb{v})\cdot \mb{s}=\mb{g}(\mb{v}_0)\cdot \mb{s}. \label{inf-attained.eqn}
 \ee
 By the oddness of $\mb{g}(\mb{v})$ in $\mb{v}$, it can be shown that the  $\mb{g}(\mb{v}_0)\cdot \mb{s}<0$ (see the related result \tb{R3} in Lemma 5.3).
 Now we have that
 \be
 -\infty~\leq ~M(\mb{s})\leq  E_P(f(y,\mb{w},\mb{0},\mb{v_0})\geq 0) +\mb{g}(\mb{v_0})\cdot \mb{s} 
 \to -\infty ~(a.s.) , ~\mbox{as $\|\mb{s}\|\to \infty$}, 
 \ee
 where the second inequality follows from the definition of infimum in $M(\mb{s})$.
 \hfill \pend
 \vs
 \vs
Let $\widehat{\mb{s}}$ be a maximizer of $M(\mb{s})$. The existence of  a $\widehat{\mb{s}}$ is guaranteed by \tb{R1} and \tb{R2}. To show the tightness of $\widehat{\mb{s}}$, it suffices to show its
measurability (see page 8 of VDV98). The latter is straightforward (see page 197 of P84, or pages 295-296 of Mass\'{e}, 2002, for example). Now we have to show that $\widehat{\mb{s}}$ is unique.
Recall that $M(\mb{s}, \mb{v})=E_P(f(y,\mb{w},\mb{0},\mb{v})\geq 0)+\mb{g}(\mb{v})\cdot\mb{s}$.
Define
$$\mathscr{V}(\widehat{\mb{s}}):=\{\mb{v}\in V(\mb{0}), M(\widehat{\mb{s}})=M(\widehat{\mb{s}},\mb{v})\}, $$
which is clearly non-empty.
Suppose that $\widehat{\mb{t}}$ is another maximizer of $M(\mb{s})$, then by \tb{R2}, $\alpha\widehat{\mb{s}}+(1-\alpha)\widehat{\mb{t}}$
 is also a maximum point for every $\alpha \in [0, 1]$.  Following Nolan, 1999, one can show that
 \vs
 \noindent
 \tb{Lemma 5.3} If \tb{A2} and \tb{A3} hold, then\vs
\noin
\tb{R3}: $\inf_{\mb{v}\in\mathscr{V}(\widehat{\mb{s}})} \mb{v}'x\leq 0$, $\forall ~x\in\R^p$;~~
\tb{R4}: $\mathscr{V}(\alpha\widehat{\mb{s}}+(1-\alpha)\widehat{\mb{t}})=\mathscr{V}(\widehat{\mb{s}})\cap \mathscr{V}(\widehat{\mb{t}})$, $\forall ~\alpha\in(0,1)$. \hfill \pend
\vs
Equipped with the results  above, we now are in the position to show that
\vs
\noindent
\tb{Lemma 5.4} If \tb{A2} and \tb{A3} hold, then $\widehat{\mb{s}}$ is unique.
\vs
\noindent
\tb{Proof}:~~
Define
$\displaystyle
\mathscr{G}:=\mbox{span}\left(\big\{\mb{g}(\mb{v}): \mb{v}\in \mathscr{V}(\alpha\widehat{\mb{s}}+(1-\alpha)\widehat{\mb{t}}, ~\mbox{for an} ~\alpha\in(0,1) \big\}\right).
$
Let $r$ be the dimension of  $\mathscr{G}$. In the sequel, consider different cases of  $r$.\vs

 If $r=1$, then there exists a $\mb{v} \in \mathscr{V}(\alpha\widehat{\mb{s}}+(1-\alpha)\widehat{\mb{t}})$ such that $\mathscr{G}$ is spanned by  $\mb{g}(\mb{v})$.  \vs

Note that $E_P(f(y,\mb{w},\mb{0},\mb{v})\geq 0)=-E_P(f(y,\mb{w},\mb{0},\mb{-v})\geq 0)$ and $\mb{g}(\mb{-v})=-\mb{g}(\mb{v})$.
  Now, by the definitions of
  $M(\mb{s})$, $M(\mb{s}, \mb{v})$ and $\mathscr {V}(\widehat{\mb{s}})$ and Lemma 5.3, we have
\bee
&&M(\widehat{\mb{s}}) = M(\widehat{\mb{s}}, \mb{v})
=\inf_{\mb{v}\in \mathcal{V(\mb{0})}} M(\widehat{\mb{s}}, \mb{v})
=-\sup_{\mb{v}\in\mathcal{V(\mb{0})}}(-M(\widehat{\mb{s}}, \mb{v}))\\[1ex]&&=-\sup_{\mb{v}\in\mathcal{V(\mb{0})}}M(\widehat{\mb{s}}, -\mb{v})
\leq -\sup_{\mb{v}\in\mathcal{V(\mb{0})}}M(\widehat{\mb{s}})=-M(\widehat{\mb{s}}),
\ene
which implies that $M(\widehat{\mb{s}})\leq 0$. By (\ref{ms.eqn}) and definitions of $M(\mb{s}, \mb{v})$ and $\widehat{\mb{s}}$, $M(\widehat{\mb{s}})\geq 0$, we conclude that $M(\widehat{\mb{s}})= 0$.
  This further implies that  $\forall~ \mb{u} \in V(\mb{0})$ and $\forall~ \mb{s} \in K$,
\be
E_P(f(y,\mb{w},\mb{0},\mb{u})\geq 0)+\mb{g}(\mb{u})\cdot \mb{s}=0 \mbox{~and~} V(\mb{0})=\mathscr{V}(\widehat{\mb{s}}).
 \label{no-name.eqn}
\ee
Now assume that there is another vector $\mb{v}_1 (\neq \pm\mb{v})\in V(\mb{0})$, then $\mb{g}(\mb{v}_1)=k\mb{g}(\mb{v})$ for some constant $k$; otherwise
$\mb{g}(\mb{v}_1)$ and $\mb{g}(\mb{v})$ are linearly independent. (\ref{no-name.eqn}) implies that
$$E_P(f(y,\mb{w},\mb{0},\mb{v}_1)\geq 0)= k E_P(f(y,\mb{w},\mb{0},\mb{v})\geq 0).$$
Write $X$ and $Y$ for $E_P(f(y,\mb{w},\mb{0},\mb{v}_1)\geq 0)$ and $E_P(f(y,\mb{w},\mb{0},\mb{v})\geq 0)$, respectively.
Then by P84 (page 149), $X$ and $Y$ have a joint bivariate normal distribution. This, however, is impossible (see (\ref{ell.eqn})) since the covariance matrix between $X$ and $Y$ has no inverse.  This implies that $\mbs^{p-1}=V(\mb{0})=\{\mb{v},-\mb{v}\}$, which can happens only if $p=1$. Namely,
both $\mb{g}(\mb{v})$ and $\widehat{\mb{s}}$ are one-dimensional. 
The uniqueness of $\widehat{\mb{s}}$  
follows in a straightforward fashion from (\ref{no-name.eqn}). 
\vs
 We now assume that $2\leq r\leq p$. Assume that $\mb{g}(\mb{v_1}),\cdots,\mb{g}(\mb{v_r})$ are linearly independent and belong to $\mathscr{G}$ and $\mb{v}_i\in \mathscr{V}(\alpha\widehat{\mb{s}}+(1-\alpha)\widehat{\mb{t}})$ for an $\alpha\in(0,1)$.
Let $S$ be any 
space that contains both $\widehat{\mb{s}}$ and $\widehat{\mb{t}}$, then both $\widehat{\mb{s}}$ and $\widehat{\mb{t}}$ satisfy the following  linear system of equations:
$$-\mb{g}(\mb{v_i})\cdot\mb{s}=E_P(f(y,\mb{w},\mb{0},\mb{v_i})\geq 0)-M(\widehat{\mb{s}}),~~ i=1,\cdots,r, ~~\mb{s}\in S$$
which immediately implies that $\widehat{\mb{s}}-\widehat{\mb{t}}=0$ is the only solution of the linear system $-\mb{g}(\mb{v_i})\cdot(\widehat{\mb{s}}-\widehat{\mb{t}})=0$, $i=1,\cdots, r$. That is, $\widehat{\mb{s}}$ is unique.  \hfill \pend
\vs
We have verified (B) completely.
As we noticed above  $\widehat{\mb{s}}_n:=n^{1/2}\bs{\beta}^*_n$  maximizes $M_n(\mb{s})$. 
To verify (A) and thus complete the proof of the theorem, we need only show that $M_n(\mb{s})\xrightarrow{d}M(\mb{s})$ uniformly in $\mb{s}\in K$, where $K\subset \R^p$ is a compact set. Note that by (\ref{taylor.eqn})
\begin{eqnarray}
M_n(\mb{s})=\inf_{\mb{v}\in\mbs^{p-1}}\Big(n^{1/2}\big(P(f(y,\mb{w},\mb{0},\mb{v})\geq0)-\alpha^*\big)+\mb{g}(\mb{v})\cdot\mb{s}+o(\|\mb{s}\|)&
\nonumber\\[1ex]
+E_n(f(y,\mb{w},\mb{0},\mb{v})\geq 0)+o_p(1)\Big),& \label{mns1.eqn}
\end{eqnarray}
Write
\begin{align}
&\lambda_n(\mb{v},\mb{s}):=n^{1/2}\big(P(f(y,\mb{w},\mb{0},\mb{v})\geq0)-\alpha^*\big)+\mb{g}(\mb{v})\cdot\mb{s}+E_n(f(y,\mb{w},\mb{0},\mb{v})\geq 0),\label{lambda.eqn}\\
&M_n^1(\mb{s}):=\inf_{\mb{v}\in\mbs^{p-1}}\lambda_n(\mb{v},\mb{s}).\label{mn1.eqn}
\end{align}
\vs
\noindent
\tb{Lemma 5.5} If \tb{A1}-\tb{A3} hold, then 
$M_n(\mb{s})\xrightarrow{d} M(\mb{s})$ uniformly over $\mb{s}\in K$.
 \vs
 \noindent
 \tb{Proof}: We employ two steps to prove the Lemma.
 \vs
(i) First, we show $\displaystyle\sup_{\mb{s}\in K}|M_n(\mb{s})- M_n^1(\mb{s})|=o_p(1)$. In  light of (\ref{mns1.eqn}) and (\ref{lambda.eqn}), we have
\bee
\sup_{\mb{s}\in K}\big|M_n(\mb{s})-M_n^1(\mb{s})\big|&=&\sup_{\mb{s}\in K}\big|\inf_{\mb{v}\in\mbs^{p-1}}\left(\lambda_n(\mb{v},\mb{s})+o(\|\mb{s}\|)+o_p(1)\right)-
\inf_{\mb{v}\in\mbs^{p-1}}\lambda_n(\mb{v},\mb{s})\big|\\[1ex]
&\leq& \sup_{\mb{s}\in K}\sup_{\mb{v}\in\mbs^{p-1}}|o(\|\mb{s}\|)+o_p(1)|=o_p(1),
\ene
where the last equality follows from two facts:
 (1) the term $o(\|\mb{s}\|)$ in (\ref{mns1.eqn}) is $o(1)$ uniformly in $\mb{s}$ over $K$, and (2) the term $o_p(1)$ in (\ref{mns1.eqn}) holds uniformly in $\mb{s}$ over $K$ for large enough $n$,  because it is obtained from application of stochastic equicontinuity over a class of functions whose members are close enough in the sense that each other is within a distance $\delta>0$  w.r.t. seminorm $\rho_P$ (see Lemma VII. 15 of P84). Thus (i) follows.
\vs
(ii) Second, we show that $M_n(\mb{s})\xrightarrow {d}M(\mb{s})$ uniformly over $\mb{s}\in K$. By virtue of (i), it suffices
to show that $M^1_n(\mb{s})\xrightarrow {d}M(\mb{s})$ uniformly over $\mb{s}\in K$. Notice that by \tb{A2},  $V(\mb{0})=\mbs^{p-1}$ and
$P(f(y,\mb{w},\mb{0},\mb{v}))-\alpha^*=0$ for any $\mb{v}\in V(0)$. Therefore,
\bee
M^1_n(\mb{s})\!\!\!&=&\!\!\!\!\inf_{v\in V(\mb{0})}\!\!\big(E_n(f(y,\mb{w},\mb{0},\mb{v})\geq 0)+\mb{g}(\mb{v})\cdot\mb{s}+n^{1/2}(P(f(y,\mb{w},\mb{0},\mb{v})\geq 0)-\alpha^*)\big)\\[1ex]
\!\!\!&=&\!\!\!\!\inf_{v\in V(\mb{0})}\!\!\big(E_n(f(y,\mb{w},\mb{0},\mb{v})\geq 0)+\mb{g}(\mb{v})\cdot\mb{s}\big)
\xrightarrow{d} \inf_{v\in V(\mb{0})}\!\!\big(E_P(f(y,\mb{w},\mb{0},\mb{v})\geq 0)+\mb{g}(\mb{v})\cdot\mb{s}\big),
\ene
where the last step follows from  the central limit theorem for empirical process (Theorem VII. 21 of P84)  and the continuous mapping theorem. The steps above hold uniformly for $\mb{s}\in K$. (A) has been verified completely. This completes the proof of the theorem 5.1.
\hfill
\pend

\vs
\noindent
\tb{Remarks 5.1}\vs
(I) Sufficient conditions for the uniqueness of $\bs{\beta^*} $ have been given at the end of Section 3. In light of Remark 4.3, a sufficient conditions for
 \tb{C2}' to hold w.r.t. RD$(\bs{\beta};P)$  have been given in Remarks 4.3. All these conditions are satisfied for the $F_{\mb{Z}}$ in Examples 2.1 and 2.2. 
\vs
(II) \tb{A2} holds true for symmetric distributions such as those regression symmetric about $\bs{\theta}$ (in this case, RD$(\bs{\theta};P)=\alpha^*$ and $V(\bs{\theta})=\mbs^{p-1}$, see Lemma 4 of RS04), which implies that the assumption \tb{A2} in the theorem could be dropped for such $F_{\mb{Z}}$.  \tb{A2} also holds for Examples 2.1 and 2.2, where $\bs{\beta^*}=(0,0)'$ and $\alpha^*=1/2$. Furthermore, if $F_{\mb{Z}}$ has a positive density and \tb{A2} holds, then $\alpha^*$=1/2, and $\bs{\beta}^*$ is unique.
\vs
(III) Theorem 5.1 could be adapted to cover the location counterpart (maximum halfspace depth estimator (aka Tukey median)), The assumptions \tb{A1}-\tb{A3} and \tb{C2}' hold
under the conditions given in Nolan, 1999 and BH99.\vs
(IV) Utilizing a different approach, BH99 treated the limit distribution of $\bs{\beta^*_n}$. BH99 skipped the verification of the two key conditions ((W1) and (W3)) in their uniqueness lemma 3.3, nevertheless.
\hfill \pend
\vs
\vspace*{-5mm}
\begin{center}
{\textbf{\large Acknowledgments}}
\end{center}
\vspace*{-2mm}
The author thanks Hanshi Zuo, Wei Shao, and Professor Emeritus James Stapleton for their careful proofreading
 and an anonymous referee and the Co-Editor-in-Chief Yimin Xiao for their insightful comments and suggestions, all of
which have led to improvements. 

\begin{center}
{\Large \tb{Appendix: Stochastic Equicontinuity and VC-class of sets}}
\end{center}
{\small The main reference of this appendix is Pollard (1984) (P84). Similar materials could also be found in Van
Der Vaart (1998) (VDV98) and Van Der Vaart and Wellner (1996) (VW96).}
\vs
\noindent
\tb{Stochastic equicontinuity}\vs
\emph{Stochastic equicontinuity} refers to a sequence of
stochastic processes $\{Z_n(t): t \in T\}$ whose shared index set $T$ comes equipped
with a semi metric $d(\cdot, \cdot)$. (a semi metric has all
the properties of a metric except that $d(s, t) = 0$ need not imply that $s$ equals
$t$.)
\vs
\noin
\tb{Definition 1} [IIV. 1, Def. 2, P84 ].  Call ${Z_n}$ stochastically equicontinuous at $t_0$  if for each $\eta > 0$
and $\epsilon > 0$ there exists a neighborhood $U$ of $t_0$ for which
\be
\limsup P\left(\sup_{U} |Z_n(t) - Z_n(t_0) | > \eta\right) < \epsilon.  \label{se.eqn}
\ee
\hfill~~~~~~~~~~~~~\pend
\vs
Because stochastic equicontinuity bounds $Z_n$ uniformly over the neighborhood
$U$, it also applies to any randomly chosen point in the neighborhood. 
If ${\tau_n}$ is a sequence of random elements of $T$ that converges in probability
to $t_0$, then
\be
Z_n(\tau)-Z(t_0)\to 0\mbox{~in probability,}
\ee
because, with probability tending to one, $\tau_n$ will belong to each $U$.
The form above will be easier to apply, especially when behavior of a particular ${\tau_n}$ sequence is under investigation.
This also is the form used in the Theorem 5.1.\vs

To establish (\ref{se.eqn}), we need the \emph{chaining technique} to prove maximum inequalities, which involves the \emph{covering number} (IIV. 2, P84).
Chaining is a technique for proving maximal inequalities for stochastic
processes, the sorts of things required if we want to check the stochastic
equicontinuity condition defined in Definition 1. It applies to any process
$\{Z(t): t\in T\}$ whose index set is equipped with a semimetric $d(\cdot,\cdot)$ that
controls the increments:
$$P\left( | Z(s) - Z(t) | > \eta\right)\leq \Delta(\eta, d(s, t)) ~\mbox{for}~ \eta > 0.$$
It works best when $\Delta(\cdot, \cdot)$ takes the form
$$\Delta(\eta, \delta) = 2 \exp( - \frac{1}{2}\eta^2/ D^2 \delta^2),$$
with $D$ a positive constant. Under some assumptions about covering numbers
for $T$, the chaining technique will lead to an economical bound on the tail
probabilities for a supremum of $\| Z(s) - Z(t)\|$ over pairs $(s, t)$.\vs
\noin
\tb{Covering number}\vs
\noin
\tb{Definition 2} [IIV. 2, Def. 8, P84].
The covering number  $N(\delta, d, T)$  is the size of the smallest $\delta$-net for $T$. That is, $N(\delta, d, T)$ equals the
smallest $m$ for which there exist points $t_1,\cdots, t_m$ with $\min_{i} d(t, t_i)\leq \delta$  for
every $t \in T$. The associated covering integral is
\be
J(\delta, d, T) = \int_{0}^{\delta} [2 \log(N(\delta, d, T)^2/u)]^ {1/2} du ~~\mbox{for~}~ 0 < \delta< 1.\label{cn-1.eqn}
\ee
\vs
\noin
\tb{Chaining Lemma} [VII. 2. Lemma 9,P84].
Let $\{Z(t): t \in T\}$ be a stochastic process whose index
set has a finite covering integral $J(\delta, d, T)$. Suppose there exists a constant $D$ such
that, for all $s$ and $t$
$$ P\big(| Z(s) - Z(t) | > \eta d(s, t)\big) \leq 2 \exp(-\eta^2/D^2)~~\mbox{for}~~ \eta > O.$$
Then there exists a countable dense subset $T^*$ of $T$ such that, for $0 < \epsilon < 1$,
$$P\Big( |Z(s) - Z(t) | > 26DJ(d(s, t))~ \mbox{for some $s$, $t$ in $T*$ with $d(s, t) \leq \epsilon$}\Big) \leq 2\epsilon $$
We can replace $T^*$ by $T$ if $Z$ has continuous sample paths. \hfill\pend
\vs
\vs
\noin
\tb{Random Covering Numbers}\vs
The symmetrization method (II. 3, P84) relates $P_n - P$ to the random
signed measure $P^o_n$ that puts mass $\pm n^{-1}$ at each of $\xi_1 ... , \xi_n$ (random sample from P), the signs being
allocated independently plus or minus, each with probability $1/2$ (see page 15 of P84). For central limit theorem calculations it is neater to work with the symmetrized empirical process $E^{o}_n =  n^{1/2} P^{o}_n$. Hoeffding's Inequality gives the clean
exponential bound for $E^{o}_n$ conditional on everything but the random signs.
For each fixed function $f\in \mathscr{F}$, a class of functions,
\begin{align}
&P\big(|E^{o}_nf| >\eta\Big|\bs{\xi}~\big) = P \Big(\Big|\sum_{i=1}^{n} \pm f(\xi_i)\Big| > \eta n^{1/2}\Big|\bs{\xi}\Big)\nonumber\\
&\leq 2 \exp \Big[ -2(\eta n^{l/2})^2\big/\sum_{i=1}^n4 f(\xi_i)^2 \Big] = 2 \exp\big[ -\frac{1}{2}\eta^2/P_nf^2\big].
\end{align}
That is, if distances between functions are measured using the $\mathscr{L}^2(P_n)$
seminorm then tail probabilities of $E^{o}_n$~ under $P(·\big|\bs{\xi})$ satisfy the exponential
bound required by the Chaining Lemma, with $D = 1$. For the purposes of
the chaining argument, $E^{o}_n$ will behave very much like the gaussian process
$B_P$\!\! , except that the bound involves the random covering number
calculated using the $\mathscr{L}^2(P_n)$ seminorm ($\rho_{P_n}(f,g)=(\int(f-g)^2dP_n)^{1/2}$, for $f,g\in \mathscr{F}$). 
Write
$$
J_2((\delta,P_n, \mathscr{F}) = \int_{0}^{\delta}\big[2\log(N_2(u, P_n,\mathscr{F})^2/u)\big]^{1/2} du
$$
for the corresponding covering integral, where we interpret P as standing for $\mathscr{L}^2(P)$ semimetrics on $\mathscr{F}$, 
the notation  $N_2(\delta, P_n, \mathscr{F})$ (a random number) agrees with Definition 2 (see II.6, Def. 32, P84).
\vs
Stochastic equicontinuity of the empirical processes $\{E_n\}$ (the signed measure $n^{l/2}(P_n - P)$) at a function
$f_0$ in $\mathscr{F}$ means roughly that, with high probability and for all $n$ large enough,
$\big| E_n f - E_n f_0 \big|$ should be uniformly small for all $f$ close enough to $f_0$. Here
closeness should be measured by the $\mathscr{L}^2(P)$ seminorm $\rho_P$. Of course we need $\mathscr{F}$ to be permissible (see Appendix C, Def. 1, P84), i.e. no measurability issue.\vs

\vs
\noin
\tb{Equicontinuity Lemma} [IIV. 4, Lemma 15, P84]\vs
Let ${\mathscr{F}}$ be a permissible class of functions with
envelope $F$ in ${\mathscr{L}}^2(P)$ (call each measurable $F$
such that $|f | \leq F$, for every $f \in \mathscr{F}$, an envelope for $\mathscr{F}$). Suppose the random covering numbers satisfy the uniformity
condition: for each $\eta > 0$ and $\epsilon > 0$ there exists a $\gamma > 0$ such that
\be
 \limsup P\big(J_2(\gamma, P_n,\mathscr{F}) > \eta \big) < \epsilon.\label{uniformity.eqn}
\ee
Then there exists a $\delta > 0$ for which
\be
\limsup P\Big(\sup_{[\delta]} |E_n(f - g)| > \eta\Big) < \epsilon,
\ee
where $[\delta] = \{(f, g): f, g \in {\mathscr{F}} ~\mbox{and}~ \rho_P(f-g)\leq \delta \}$. \hfill\pend

\vs
Up to this point, there are two approaches to establish the stochastic equicontinuity: (i) via Definition 1 (ii) via Equicontinuity Lemma.
The first approach is usually more challenging, the second one is equivalently to verify the  uniformity
condition for the random covering numbers. \vs

A sufficient condition for the latter is the graphs of the functions in $\mathscr{F}$ have polynomial discrimination.
The graph of a real-valued function $f$ on a set $S$ is defined as the subset
$$G_f = \{(s, t): 0\leq t \leq f(s) ~\mbox{or}~ f(s)\leq t \leq 0, s \in S \}.$$
\vs
If the graphs of the functions in $\mathscr{F}$ have polynomial discrimination, then $N_2(u, P_n, \mathscr{F})$ is bounded by a polynomial $A(u^{-1})^W$ in $u^{-1}$ with $A$ and $W$ not depending on
$P_n$ (Lemma II. 36, P84), which amply suffices for the Equicontinuity Lemma: For each $\eta>0$, there is a $\gamma>0$ so that
$J_2(\gamma, P_n, \mathscr{F})\leq \eta$ for every $P_n$.  Therefore, the graphs of the functions in $\mathscr{F}$ having polynomial discrimination becomes a key point for Equicontinuity Lemma.
\vs
\vs
\noin
\tb{Polynomial discrimination}\vs
\noin
\tb{Definition 3} [II.4, Def.13, P84]. Let $\mathscr{D}$ be a class of subsets of some space $S$. It is said to have
polynomial discrimination (of degree $v$) if there exists a polynomial $\rho(\cdot)$
(of degree $v$) such that, from every set of $N$ points in $S$, the class picks out at
most $\rho(N)$ distinct subsets. Formally, if $S_0$ consists of $N$ points, then there
are at most $\rho(N)$ distinct sets of the form $S_0\cap D$ with $D \in {\cal{D}}$. Call $\rho(\cdot)$ the
discriminating polynomial for ${\mathscr{D}}$. ${\mathscr{D}}$ is also called  a VC-class of sets (see Vapnik and Chervonenkis, 1971).
\vs
\noin
\tb{Generalized  Glivenko-Cantelli theorem}\vs
\noin
\tb{Theorem 1} [II.4, Th.14, P84]. Let $P$ be a probability measure on a space $S$. For every permissible
class ${\mathscr{D}}$ of subsets of $S$ with polynomial discrimination,
$$\sup_{{\mathscr{D}}}|P_nD - PD|\to 0~~ \mbox{almost surely}.
$$\vs
\noindent
\tb{Examples}\vs
\tb{1}. Let ${\mathscr{D}}=\{(-\infty, t], t\in \R \}$. The collection of sets is the one in the traditional Glivenko-Cantelli theorem in one-dimension. The ${\mathscr{D}}$ can pick at most $(n+1)$ subsets for any set of $n$ points in the line.  $\mathscr{D}$ has polynomial discrimination. Theorem 1 holds true.
\vs
\tb{2}. Let $\mathscr{D}=\{(-\infty, \bs{t}], \bs{t}\in \R^2 \}$. The collection of all quadrants of the form $ (-\infty, \bs{t}]$ in $\R^2$, which can pick fewer than $(n+1)^2$ subsets from a set of $n$ points in the plane. $\mathscr{D}$ has polynomial discrimination. Theorem 1 holds true.
\vs
\tb{3}. Let $\mathscr{D}$ be the class of all closed halfspaces in $\R^d$, then it can pick at most $O(n^2)$ subsets from a set of $n$ points in $\R^d$.
 $\mathscr{D}$ has polynomial discrimination. Theorem 1 holds true.\vs
\tb{4}. Let $\mathscr{D}$ be the class of closed, convex sets. From every collection
of n points lying on the circumference of a circle in $\R^2$, it can pick out all $2^n$ subsets. $\mathscr{D}$ no longer has polynomial discrimination. \hfill\pend

\vs

Back to the Equicontinuity Lemma, a sufficient condition for the uniformity condition in the lemma is that the graphs of the functions in ${\mathscr{F}}$ have polynomial discrimination. How to verify the latter becomes the key point. It turns out that this can be done straightforwardly by the following lemma.\vs
\noindent
\tb{Lemma 1} [II.5, Lemma 28, P84]. Let ${\mathscr{F}}$ be a finite-dimensional vector space of real functions on S.
The class of graphs of functions in ${\mathscr{F}}$ has polynomial discrimination.
\vs
The following lemmas are equally useful as well.
\vs
\noindent
\tb{Lemma 2} [II.4, Lemma 18, P84]. Let ${\mathscr{G}}$ be a finite-dimensional vector space of real functions on $S$.
The class of sets of the form $\{g \geq 0\}$, for g in ${\mathscr{G}}$, has polynomial discrimination of degree no greater than the dimension of ${\mathscr{G}}$.
\vs
\noindent
\tb{Lemma 3} [II.4, Lemma 15, P84].  If  ${\mathscr{C}}$ and ${\mathscr{D}}$ have polynomial discrimination, then so do each of:
(i) $\{D^c: D\in {\mathscr{D}}\}$; (ii) $\{C \cup D: C\in {\mathscr{C}} ~ \mbox{and~} D \in  {\mathscr{D}}\}$;
(iii) $\{C \cap D: C\in {\mathscr{C}} ~ \mbox{and~} D \in  {\mathscr{D}}\}$.
\vs
\vs
\noindent
\tb{Example}\vs
\tb{5}. Consider three classes of functions that appear in \tb{Theorems 5.1} and \tb{4.2}, 
$${\mathscr{F}}_{\mb{v}}=\{f(\cdot, \cdot, \mb{0}, \mb{v}), \mb{v}\in \mbs^{p-1}\},$$
$$\mathscr{F}_{\bs{\beta}}=\{f(\cdot, \cdot, \bs{\beta}, \mb{v}_0), \bs{\beta}\in \R^p, \mb{v}_0\in \mbs^{p-1}~\mbox{is fixed}\},$$
$$\mathscr{F}_{\bs{\beta}, \mb{v}}=\{f(\cdot, \cdot, \bs{\beta}, \mb{v}), \bs{\beta}\in \R^p, \mb{v}\in \mbs^{p-1}\}.$$
Under \tb{A0} given in Theorem 4.2, the three classes have a square integrable envelope F. 
We want to show that the graphs of the functions in these classes have polynomial discrimination.
\vs
It suffice to show this for $\mathscr{F}_{\bs{\beta}, \mb{v}}$ since other are just special cases.
\vs
The graph of a function in $\mathscr{F}_{\bs{\beta}, \mb{v}}$ contains a point $((y,\mb{w}), t)$  if and only if
$0\leq f(y,\mb{w}, \bs{\beta}, \mb{v})\leq t$ or $t\leq f(y,\mb{w}, \bs{\beta}, \mb{v})\leq 0$.
Therefore, the total number of subsets of given $n$ points $((y_i, \mb{w}_i), t_i)$ that can be picked out by the graphs of functions in $\mathscr{F}_{\bs{\beta}, \mb{v}}$ is less than the total number of those picked out by the union of classes of sets $\{f(y_i, \mb{w}_i, \bs{\beta}, \mb{v})\geq 0\}\cup \{f(y_i, \mb{w}_i, \bs{\beta}, \mb{v})\leq0\}$. We now show that each class of sets has polynomial discrimination, so does the union (by Lemma 3) and consequently so do the graphs of functions in $\mathscr{F}_{\bs{\beta}, \mb{v}}$.\vs
Note that $f(y,\mb{w}, \bs{\beta}, \mb{v})=(y-\bs{\beta}'\mb{w})\mb{v}'\mb{w}$ with $\mb{w}'=(1, \mb{x}')$. Hence
\bee
&&\{f(y,\mb{w}, \bs{\beta}, \mb{v})\geq 0, \bs{\beta}\in \R^p, \mb{v}\in \mbs^{p-1} \}=
\Big\{\Big(\{(y-\bs{\beta}' \mb{w})\geq0\}\cap\{\mb{u}' \mb{x}< v\}\Big)\\ &&\cup \Big(\{ (y-\bs{\beta}' \mb{w})<0\}\cap\{\mb{u}' \mb{x}\geq v\}\Big),   \mb{u}\in \mbs^{p-2}, v\in \R^1, \bs{\beta}\in \R^p\Big\}.
\ene
The RHS is built up from sets of the form $\{g\geq 0\}$ with $g$ in the finite-dimensional
vector space of functions. There are four classes of functions on the RHS each forming a finite dimensional vector space of real functions.
By Lemmas 1 and 3, we conclude that the graphs of the functions in $\mathscr{F}_{\bs{\beta}, \mb{v}}$ have polynomial discrimination. So do the other class of functions.
\vs
 That is, they are VC class, 
or polynomial class of functions in the terminology of P84. 
We have all needed for invoking Equicontinuity Lemma. 
However, to invoke the Corollary 3.2 in KP90, as did in the proof of Theorem 4.2, we actually need to show that  $\mathscr{F}_{\bs{\beta}, \mb{v}}$ is a
manageable class of functions, a notion defined in Pollard, 1989 (P89).  On the other hand, every VC class is a Euclidean class (a notion introduced in Nolan and Pollard, 1987), luckily enough, every Euclidean class is manageable (P89). \hfill\pend

\end{document}